\documentstyle[11pt]{amsart}
\def\R{{\Bbb R}}

\def\mm{{M}}
\def\ii{{I}}
\def\jj{{J}}
\def\kk{{K}}
\def\gg{{G}}
\def\uu{{U}}

\def\ww{{W}}
\def\hh{{H}}
\def\h{{h}}
\def\ccj{{C^1}}
\def\ep{{\varepsilon}}
\def\xbar{{\bar{\x}}}
\def\fita{{\varphi_\ta}}
\def\ddom{\mathop{\rm Dom}\nolimits}
\def\ccodom{\mathop{\rm Codom}\nolimits}
\def\ggal{{\gg_{\al}}}

\def\gg-al{{\gg_{-\al}}}

\def\ggal-be{{\gg_{\al-\be}}}
\def\ggta{{\gg_{\ta}}}
\def\ggta-ro{{\gg_{\ta-\ro}}}

\def\Rn{{{\Bbb R}^n}}

\def\ro{{\rho}}

\def\si{{\sigma}}

\def\ta{{\tau}}
\def\al{{\alpha}}
\def\be{{\beta}}
\def\ga{{\gamma}}

\def\om{{\omega}}

\def\ep{{\varepsilon}}

\def\xdot{{\dot{x}}}
\def\x{{x}}

\def\gg{{G}}
\def\h{{h}}
\def\hh{{H}}
\def\f{{f}}
\def\ff{{F}}

\def\kk{{K}}

\def\firo{{\varphi_\ro}}

\def\fftasi{{F_{\ta\si}}}
\def\fftaro{{F_{\ta\ro}}}
\def\ffsiro{{F_{\si\ro}}}
\def\ffsita{{F_{\si\ta}}}
\def\ffsiom{{F_{\si\om}}}
\def\ffrosi{{F_{\ro\si}}}
\def\ffroro{{F_{\ro\ro}}}
\def\ffomom{{F_{\om\om}}}
\def\ffomsi{{F_{\om\si}}}
\def\ffomro{{F_{\om\ro}}}

\def\ffsisi{{F_{\si\si}}}

\begin{document}
\sloppy

\pagestyle{plain}
\setcounter{page}{1}
\newtheorem*{veta}{Theorem} 
\newtheorem*{lemma}{Lemma}
\frenchspacing

\title{The Cauchy atlas on the manifold of all complete ODE solutions}
\author{Petr Chl\'{a}dek, Lubom\'\i r Klapka}
\curraddr{Mathematical Institute of the Silesian University at Opava\\
\newline
\indent
Na Rybn\'\i\v{c}ku 1\\
746 01 Opava\\Czech Republic}
\makeatletter
\email{Petr.Chladek@math.slu.cz, Lubomir.Klapka@math.slu.cz}
\makeatother
\begin{abstract}
\noindent
In this paper the necessary and sufficient conditions for a mapping to be the dependence of the complete
solution of some
$\ccj$ first-order ordinary differential equation on the initial Cauchy condition
are deduced. The result is obtained by studying the Cauchy atlas on a manifold
of complete solutions. The proof is constructive - the corresponding differential equation is obtained.
The autonomous case and the linear case are discussed. The relation to the
Sincov functional equation is clarified.
\\


\\
\end{abstract}

\maketitle

\vspace{-1cm}
\section{Introduction}
\noindent
In Section 2, we introduce the first-order ordinary differential equation and the Cauchy atlas on a manifold of
complete solutions. Four necessary conditions for a mapping to be the coordinates transformation of
the Cauchy atlas are given.

In Section 3, we also prove the sufficient conditions (\ref{3}) for a system $\ff=\{ \fftasi \}_{\ta , \si \in \R}$
of functions to be the Cauchy atlas coordinates transformations.
The necessary conditions are well-known and included
in many papers and monographs (see, e.g., \cite{arnold}, \cite{pontrjagin}), whereas the sufficient conditions
are presented here, as far as the authors know, for the first time. These conditions are in the form of the
composition of iterations in the non-autonomous dynamical systems theory (see, e.g., \cite{snoha}). 

Sections 4, 5 and 6 contain special cases of the first-order ordinary differential equations. More precisely,
we discuss here autonomous equations and linear inhomogeneous equations.
Section 7 contains a simple example.
\bigskip

\section{The Cauchy atlas}
\noindent
Let $\mm$ be a set, $n$ a natural number, $\ii$ an index set. A {\it coordinates chart} is a
bijection; its domain is a subset of $\mm$ and its codomain is
a subset of $\Rn$. A system of coordinates charts 
$\{ \fita \}_{\ta \in I}$ is an {\it atlas} if and only if the
union of domains of all coordinates charts is equal to $\mm$. 
A {\it manifold} is a set endowed with an atlas. 
By {\it coordinates transformation} on a manifold we mean the mapping
\begin{equation}
\begin{array}{c}
\label{12}
\fftaro \colon \firo(\mbox{Dom}(\fita) \cap \mbox{Dom}(\firo) ) \ni a \mapsto \qquad\qquad\qquad\qquad\\
\qquad\qquad\qquad \fita(\firo^{-1}(a)) \in \fita(\mbox{Dom}(\fita) \cap \mbox{Dom}(\firo) ),
\end{array}
\end{equation}
where $\fita$ and $\firo$ are coordinates charts. In the definition of the coordinates transformation we admit
the empty mapping, i.e. $\mbox{Dom}(\fita) \cap \mbox{Dom}(\firo)=\emptyset$.

Let $\f$ be a $\ccj$-mapping, where $\mbox{Dom}(\f)$ is a non-empty open subset of $\R \times \Rn$ and 
$\mbox{Codom}(\f)=\Rn$.
Consider the first-order ordinary differential equation
\begin{equation}
\label{1}
\xdot=\f(\ta,\x)
\end{equation}
and the Cauchy initial condition
\begin{equation}
\label{2}
\x(\ro)=a,
\end{equation}
where $(\ro,a) \in \mbox{Dom}(\f)$.

Let $\mm$ be the set of all complete solutions
of $\ccj$ equation (\ref{1}).
For each $\ta \in \R$ we define mapping $\fita$, where
$\mbox{Dom}(\fita)=\{\x \in \mm \mid \ta \in \mbox{Dom}(\x) \}, 
\mbox{Codom}(\fita)=\{ a \in \Rn \mid (\ta,a) \in \mbox{Dom}(\f) \}$ and
$\fita$ maps a complete solution $\x$ to the value $\x(\ta)$.
From the existence and the uniqueness of the complete solution of the equation (\ref{1}), (\ref{2})
(see, e.g., \cite[chapter V]{hartman}) the bijectivity of $\fita$ follows.
Hence $\fita$ is a coordinates chart on $\mm$.
The system $\{\fita\}_{\ta \in \R}$ forms a {\it Cauchy atlas} on $\mm$.

\begin{lemma}
Let $\ff=\{ \fftaro \}_{\ta , \ro \in \R}$ be the coordinates transformations system of the Cauchy atlas on the set
$\mm$ of all complete solutions of $\ccj$ equation $(\ref{1})$. Then

$1.$ $\kk=\{(\ta, \ro, a) \in \R \times \R \times \Rn \mid a \in \ddom(\fftaro) \}$ is a non-empty

\quad open set,

$2.$ each function from $\ff$ is bijection,

$3.$ mapping $\kk \ni (\ta,\ro,a) \mapsto \partial \fftaro(a) / \partial \ta \in \Rn$ is $\ccj$,

$4.$ $\jj(\ro,a)=\{ \ta \in \R \mid a \in \ddom(\fftaro) \}$ is an open interval for each

\quad $\ro \in \R$, $a \in \Rn$.
\end{lemma}

\newcounter{k1} \setcounter{k1}{4}
\noindent
{\it Proof.} From (\ref{12}) the mapping 
\begin{equation}
\label{11}
\x \colon \jj(\ro,a) \ni \ta \rightarrow \fftaro(a) \in \Rn
\end{equation}
is the complete solution of (\ref{1}) satisfying (\ref{2}). The
dependence of the solution on the Cauchy condition has an open domain (see \cite[Chapter 4,
$\fnsymbol{k1} 23$]{pontrjagin}). The statement 1 follows from this and from the non-emptiness of $\ddom(\f)$. The
statement 2 is obtained from the definition of coordinates transformations. It follows from \cite[Chapter 4,
$\fnsymbol{k1} 32$]{arnold} that the mapping
$\kk \ni (\ta,\ro,a) \mapsto \fftaro(a) \in \Rn$ is $\ccj$. The mapping $\f$ is $\ccj$ by assumptions.
Therefore $(\ta,\ro,a) \mapsto f(\ta,\fftaro(a))$ is $\ccj$.
The statement 3 follows from (\ref{1}), (\ref{11}). Since $\jj(\ro,a)$ is the domain of the complete solution
(\ref{11}), we get the statement 4 from the openness of $\ddom(\f)$. $\Box$
\bigskip

\section{Transformations}
\noindent
Let $\ff=\{ \fftaro \}_{\ta , \ro \in \R}$ be a system of functions, where domain and codomain of each function
from $\ff$ are subsets of $\Rn$. The following Theorem gives the necessary and sufficient conditions
for $\ff$ to be a system of the Cauchy atlas transformations.

\begin{veta}
The functions from the system $\ff=\{ \fftasi \}_{\ta , \si \in \R}$ are coordinates transformations of the Cauchy
atlas  on a manifold  $\mm$ of the complete solutions of some $\ccj$ equation 
$(\ref{1})$ if and only if the statements $1$, $2$, $3$, $4$ from Lemma are satisfied and the condition 
\begin{equation}
\label{3}
\fftasi ( \ffsiro (a))= \fftaro (a) 
\end{equation}
holds for each $\ta,\si,\ro \in \R$ and
for each 
$a \in \ffsiro^{-1} (\ccodom (
\ffsiro) \cap \ccodom( \ffsita ) )$.
\end{veta}

\noindent
{\it Proof.} Assume that $\ff$ is a system of coordinates transformations of the Cauchy atlas on $\mm$. According
to Lemma, the statements 1, 2, 3, 4 are satisfied. From the definition of the coordinates transformation
we have the condition (\ref{3}) immediately.

Let us suppose the condition (\ref{3}) and the statements 1, 2, 3, 4 from Lemma hold.
From (\ref{3}) we get $\ffsiro(\ffrosi(a))=\ffsisi(a)$ 
for each $a \in  
\mbox{Dom}(\ffrosi)$.
If we put $\ro=\si$, then from the statement 2 of Lemma we obtain
$a=\ffsisi(a)$
for each $a \in \mbox{Dom} (\ffsisi)$.
From above
\begin{equation}
\label{8}
\mbox{Dom}(\ffrosi) \subseteq \mbox{Dom}(\ffsisi), 
\end{equation}
\begin{equation}
\label{9}
\ffrosi=\ffsiro^{-1}
\end{equation}
for each $\ro,\si \in \R$.

Consider the equation (\ref{1}), where
\begin{equation}
\label{4}
\f(\ta,a)=\left.\frac{\partial \fftasi(a)}{\partial \ta} \right |_{\si=\ta}.
\end{equation}
From the statement 3 of Lemma we see that $\f$ is $\ccj$. 
Let $(\ro,a) \in \ddom(\f)$ be fixed. Then $a \in \mbox{Dom}(\ffroro)$.
Moreover,  $\jj(\ro,a)$ is non-empty according to the statement 4 of Lemma. Let $\si \in \jj(\ro,a)$.
Then  $a \in \mbox{Dom}(\ffsiro)$. Let $\ta \in \jj(\si,\ffsiro(a))$. This set is non-empty, since from
(\ref{8}), (\ref{9}) we have $\si \in \jj(\si,\ffsiro(a))$. From (\ref{9}) $a \in \ffsiro^{-1} (\ccodom (
\ffsiro) \cap \ccodom( \ffsita ) )$. Therefore the condition (\ref{3}) holds for such $a$.
Differentiating (\ref{3}) with respect to $\ta$, putting $\si=\ta$ and using (\ref{4}) we get
$$
\frac{\partial \fftaro(a)}{\partial \ta}=\f(\ta,\fftaro(a)).
$$
Then the mapping $\x$ defined by (\ref{11}) is the solution of the equation 
(\ref{1}), (\ref{4}).
We will check that $\x$ is a complete solution. Let us suppose $\x$ is not complete. Then there
exists a complete solution $\xbar$ such that 
$\x=\left. \xbar \right |_{\jj(\ro,a)}$. 
At least one of the values $\mbox{sup}(\jj(\ro,a)), \mbox{inf}(\jj(\ro,a))$ is an element of $\mbox{Dom}(\xbar)$.
Let us suppose $\om=\mbox{sup}(\jj(\ro,a)) \in \mbox{Dom}(\xbar)$ (the case $\mbox{inf}(\jj(\ro,a)) \in
\mbox{Dom}(\xbar)$ is analogous).
Then
$(\om,\xbar(\om)) \in \mbox{Dom}(\f)$. Further from (\ref{4}) $\xbar(\om) \in \mbox{Dom}(\ffomom)$. 
Therefore $(\om, \om, \xbar(\om)) \in \kk$.
By the statement 1 of Lemma, there
exists $\ep>0$ such that for each $\si$ from an open interval $(\om - \ep, \om)$ we have 
$(\om,\si,\xbar(\si)) = (\om,\si,\x(\si)) = (\om, \si, \ffsiro(a))
\in \kk$.
By using (\ref{9}) we get $a \in \ffsiro^{-1} (\mbox{Codom} (
\ffsiro) \cap \mbox{Codom}( \ffsiom ) )$.
From (\ref{3}) we obtain $\ffomsi(\ffsiro(a))=\ffomro (a)$. Therefore $\mbox{sup}(\jj(\ro,a))=\om \in
\jj(\ro,a)$. Nevertheless from the statement 4 of Lemma $\jj(\ro,a)$ is an open set. This
contradiction proves that
$\x$ is the complete solution. 
Thus for any Cauchy condition (\ref{2}), where $(\ro,a) \in \mbox{Dom}(\f)$, we have the complete solution
(\ref{11}). Since $\f$ is $\ccj$, any other complete solution does not exist.
Then we can construct the Cauchy atlas. From (\ref{11}), (\ref{9}) and from the statement 4 of Lemma the condition
(\ref{12}) holds. $\Box$
\bigskip

\section{The autonomous case}
\noindent
Let $\uu \subset \Rn$ be an open set.
Let (\ref{1}) be a $\ccj$ equation with $\ddom(\f)=\R \times \uu$, where $\f \colon (\ta,a) \mapsto \xi(a)$
and $\xi \colon \uu \rightarrow \Rn$. Then the mapping $\jj(\ro,a) \ni \ta \mapsto \ff_{\ta-\ro,0}(a) \in \Rn$
is the complete solution of this equation satisfying the same Cauchy condition as the solution (\ref{11}).
Hence we can put 
\begin{equation}
\label{a1}
\fftaro=\ggta-ro,
\end{equation}
where $\gg_{\ta}=\ff_{\ta 0}$. From (\ref{3}), (\ref{a1}) we have
$$
\gg_{\al}(\gg_{\be}(a))=\gg_{\al+\be}(a)
$$
for each $a \in \gg_{\be}^{-1} (\ccodom(\gg_{\be}) \cap \ccodom(\gg_{-\al}))$.
The mapping $\gg \colon (\al,a) \mapsto \gg_{\al}(a)$, where 
$\ddom(\gg)=\{  (\al,a) \in \R \times \Rn \mid a \in \ddom(\gg_{\al}) \}$,
is the {\it maximal flow} of the vector field $\xi$ (see, e.g., \cite[Chapter 17]{lee}).
The mappings $\gg_{\ta}$ form a {\it local one-parameter group} of transformations (for $C^{\infty}$ case see,
e.g., \cite[Section 1.2]{olver}). If $\ddom(\gg_{\ta})=\uu$ for each $\ta \in \R$, then $\gg_{\ta}$'s form a 
{\it group of  transformations} of $\uu$.
\bigskip

\section{The linear case}
\noindent
Let $\ii$ be an open interval. Let us consider the linear inhomogeneous functions $\fftasi$, where
$$
\ddom(\fftasi)=\left\{
\begin{array}{ll}
\Rn \quad & \mbox{for } \ta,\si \in \ii, \\
\emptyset \quad & \mbox{otherwise.}
\end{array}
\right.
$$ 
 The condition (\ref{3}) can be rewritten as {\it Sincov's functional equation} (see \cite[section 8.1]{aczel})
$$
\fftasi \circ \ffsiro = \fftaro
$$
for each $\ta$, $\si$, $\ro \in \ii$.
Its general solution has the following form
$$
\fftasi(a)=\ww_{\ta}(\ww_{\si}^{-1}(a)+\h_{\ta}-\h_{\si}),
$$
where $\ww_{\ta} \colon \Rn \rightarrow \Rn$ is an arbitrary linear automorphism
and $\h_{\ta}$ is an arbitrary term of $\Rn$ for each $\ta \in
\ii$. If the conditions from Theorem are satisfied, then the differential equation (\ref{1}), (\ref{4}) is
also linear inhomogeneous and $\ddom(\f)=\ii \times \Rn$.
Moreover, $\ta \mapsto \ww_{\ta}$ is the {\it Wronski matrix} and $\ta \mapsto \ww_{\ta} \h_{\ta}$ is the
{\it particular solution} of this equation.
\bigskip

\section{The linear autonomous case}
\noindent
Let the mapping $\gg_{\ta} \colon \Rn \rightarrow \Rn$ defined by (\ref{a1}) be linear inhomogeneous for each
$\ta \in \R$. We can rewrite the condition (\ref{3}) as 
\begin{equation}
\label{la1}
\gg_{\al} \circ \gg_{\be}=\gg_{\al+\be}. 
\end{equation}
Therefore $\gg_{\ta}$'s form
the group of affine transformations of $\Rn$.
Let us suppose the mapping $\be \mapsto \gg_{\be}$ is continuous. 
We define a mapping 
$$
\hh \colon \ep \mapsto \frac{1}{2\ep} \int\limits_{-\ep}^{\ep} \gg_{\be} \,d\be.
$$
Since $\lim_{\ep \rightarrow 0} \hh_{\ep}=\mbox{id}_{\Rn}$, from continuity,
there exists $\ep>0$ such that $\hh_{\ep}$ is invertible.
By integrating (\ref{la1}) and substituting $\ga=\al+\be$ we obtain
$$
\gg_{\al}=\frac{1}{2\ep} \int\limits_{\al-\ep}^{\al+\ep} \gg_{\ga} \circ \hh_{\ep}^{-1} \,d\ga.
$$
From this and from (\ref{a1}) we have the statement 3 of Lemma.
From (\ref{la1}) and from Theorem we can see that
functions from the system $\ff=\{ \gg_{\ta-\si} \}_{\ta , \si \in \R}$ are the coordinates transformations of the
Cauchy atlas  on a manifold  $\mm$ of the complete solutions of some $\ccj$ equation 
$(\ref{1})$.
From (\ref{4}), (\ref{a1}) the equation (\ref{1}) is linear inhomogeneous with constant coefficients.
\bigskip

\section{Example}
\noindent
{\bf Example.} Let us consider the equation (\ref{1}), where $\f \colon \R^2 \ni (\ta,a) \mapsto a^2 \in \R$. It is
easy to see that for this equation we have 
$$
\fftasi(a)=\frac{a}{1+(\si-\ta)a},
$$
where $\ddom(\fftasi)=\{ a \in \R \mid (\ta-\si)a<1 \}$, $\ccodom(\fftasi)=\{ a \in \R \mid (\ta-\si)a>- 1 \}$.
The mapping $\gg \colon (\ta,a) \mapsto a / (1-\ta a)$ is the maximal flow of
$\xi \colon a \mapsto a^2$.
\bigskip

\noindent
{\bf Acknowledgment.} The first author was supported by the grant MSM 192400002 from the Czech Ministry of
Education.

\renewcommand{\refname}{References}


\begin{thebibliography}{20}
\bibitem{aczel} J. Acz\'{e}l,
{\it Lectures on functional equations and their applications},
New York and London: Academic Press, 1966.
\bibitem{arnold} V.I. Arnold,
{\it Ordinary differential equations} (3rd edition) (Russian),
Moscow: Nauka, 1984.
\bibitem{hartman} P. Hartman,
{\it Ordinary differential equations},
New York-London-Sydney: John Wiley and Sons, 1964.
\bibitem{snoha} S. Kolyada, L. Snoha,
{\it Topological entropy of nonautonomous dynamical systems},
Random \& Computational Dynamics, 4(2\&3) (1996), 205-233, ISSN 1061-835X. 
\bibitem{lee} J.M. Lee,
{\it Introduction to smooth manifolds},
New York: Springer, 2002.
\bibitem{olver} P. Olver,
{\it Applications of Lie groups to differential equations} (2nd edition),
New York: Springer-Verlag, 1993.
\bibitem{pontrjagin} L.S. Pontrjagin,
{\it Ordinary differential equations} (Russian),
Moscow: Nauka, 1965.
\\
\\
\end{thebibliography}
\end{document}